\documentclass{article}
\usepackage[margin=15mm]{geometry}
\usepackage{amsmath,amssymb,amsthm,ascmac,url,bm,physics}
\usepackage{graphicx} 
\usepackage{etoolbox}
\usepackage{tikz}
\usetikzlibrary{graphs}

\usepackage{algorithm}
\usepackage{algpseudocode}

\title{The exact value of $c_1(K_{2,n})$}
\author{Hiroaki Mori}
\date{}

\newtheorem{thm}{Theorem}
\newtheorem{prop}[thm]{Proposition}

\newcommand{\N}{\mathbb{N}}
\newcommand{\Z}{\mathbb{Z}}
\begin{document}
\maketitle
\begin{abstract}
    For a graph $G$, let $c_1(G)$ be the largest distortion necessary to embed any shortest-path metric on $G$ into $\ell_1$, and for any natural number $n,m\in\N$, denote $K_{n,m}$ as the complete bipartite graph. In this note, we caculate the value of $c_1(K_{2,n})$, more precisely we prove $c_1(K_{2,n})=\frac{3k-2}{2k-1}$ where $k=\lceil\frac{n}{2}\rceil$.
\end{abstract}

\section{Introduction}
For two metric spaces $(X,d_X)$ and $(Y,d_Y)$, the \textit{distortion} of $f:X\to Y$ is defined as
\[\mathrm{dist}(f)=\left(\sup_{d_X(x,y)\neq 0}\frac{d_Y(f(x),f(y))}{d_X(x,y)}\right)\cdot\left(\sup_{d_Y(f(x),f(y))\neq 0}\frac{d_X(x,y)}{d_Y(f(x),f(y))}\right)\]

For a graph $G$ and $w:E(G)\to[0,+\infty)$, let $c_1(G,w)$ be the minimum distortion required to embed the shortest path metric of $(G,w)$ into $\ell_1$-space and define $c_1(G)=\sup_{w\text{ is weight on }G}c_1(G,w)$. This quantity is closely related to the flow-cut gap \cite{Gupta2004}.

\begin{thm}[Gutpa et al. \cite{Gupta2004}]
    The maximum ratio between the sparsest cut and the maximum concurrent flow ranges over all instances on $G$ matches $c_1(G)$.
\end{thm}

In this note, we calculate the exact value of $c_1(K_{2,n})$ for all positive numbers $n$. It strengthens the result of \cite{ostrovsky2007low} which shows $\lim_{n\to\infty}c_1(K_{2,n})=3/2$. To our best knowledge, it is only known that $c_1(K_{2,1})=c_1(K_{2,2})=1$ and $c_1(K_{2,3})=4/3$, thus this is the first example where $c_1(G)$ of an infinite number of graph families could be precisely computed, excluding the case $c_1(G)=1$.

\begin{thm}
    For any positive number $n$, we have $c_1(K_{2,n})=\frac{3\lceil\frac{n}{2}\rceil-2}{2\lceil\frac{n}{2}\rceil-1}$.
\end{thm}

On the following discussion, we fix a natural number $n\in\N$ arbitrary and set $k=\lceil\frac{n}{2}\rceil$. We also use $[n]=\{1,2,\ldots,n\}=\{k\in\mathbb{N}:1\leq k\leq n\}$, $\binom{[n]}{m}=\{A\subset [n]:|A|=m\}$ for any natural numbers $n,m\in\mathbb{N}$.

\section{Lower Bounds}
Clearly, $c_1(K_{2,n})$ is monotonically increasing with $n$. So we only have to prove $c_1(K_{2,2k+1})\geq \frac{3k+1}{2k+1}$ for any $k\in\N$. We see that embedding $K_{2,2k+1}$ equipped with unit weight needs at least $\frac{3k+1}{2k+1}$-distortion. To prove this, we employ the following inequality. See Deza and Laurent\cite{Deza} for the proof of $\ell_1$-metric is hypermetric.
\begin{thm}[Hypermetric inequality\cite{Deza}]
    Let $X\subset\ell_1$ be a finite metric space and $d$ be its metric. Then, for any sequence of integers $\{b_x\}_{x\in X}\subset\Z$ with $\sum_{x\in X}b_x=1$, we have
    \[\sum_{x,y\in X}b_xb_yd(x,y)\leq 0.\]
\end{thm}

Write $V(K_{2,2k+1})=A\cup B$ as the ordinary partition of bipartite graph, assume $|A|=2,|B|=2k+1$. Take any non-contracting embedding $f:K_{2,n}\to\ell_1$ and let $D$ be the distortion of $f$. Apply the hypermetric inequality as $b_x=-k$ for $x\in A$ and $b_x=1$ for $x\in B$, then we have
\[k^2\sum_{x,y\in A}\norm{f(x)-f(y)}_1+\sum_{x,y\in B}\norm{f(x)-f(y)}_1\leq k\sum_{x\in A,y\in B}\norm{f(x)-f(y)}_1\]
However, we have
\[k^2\sum_{x,y\in A}\norm{f(x)-f(y)}_1+\sum_{x,y\in B}\norm{f(x)-f(y)}_1\geq 2k^2+2\binom{2k+1}{2}=6k^2+2k\]
while
\[k\sum_{x\in A,y\in B}\norm{f(x)-f(y)}_1\leq 2k(2k+1)D\]
Thus, we derive $D\geq\frac{3k+1}{2k+1}$.

\section{Upper Bounds}
By the monotonicity, we only have to prove $c_1(K_{2,2k})\leq\frac{3k-2}{2k-1}$ for all $k\in\N$. Let $w:E(K_{2,2k})\to\mathbb{R}_{\geq 0}$ be an arbitrary weight of $K_{2,2k}$. We first note the following compactness argument. See Deza and Laurent\cite{Deza} for the proof.
\begin{prop}
    For any compact metric space $(X,d)$ and metric space $Y$, if $X$ can embed into $Y$ with distortion $D+\varepsilon$ for all $\varepsilon>0$, then $X$ can embed into $Y$ with distortion $D$.
\end{prop}
Therefore, by scaling and approximating $w$ with rational valued weight, we only have to deal with the graph $G$, which is obtained from subdividing $K_{2,2k}$ several times. We point out that cutting loose paths does not affect the distortion. We use the statement, as seen in Chakrabarti et al. \cite{Chakrabarti_et_al_2008}, in the following manner.
\begin{prop}\label{prop:shrink}
    Let $G$ be any bipartite graph and $d_G$ be its shortest path metric. For any $x,y\in V(G)$ and the path $P$ between $x,y$ with its length is strictly larger than $d_G(x,y)$ and all the inner points of $P$ are degree $2$.

    Then, there is a non-trival cut $R$ of inner points of $P$, which has following properties: let $H$ be the graph obtained from $G$ by shrinking the edges between $R$ and $V(G)\backslash R$, then for any embedding $f:V(H)\to\ell_1$, there exists an embedding $g:V(G)\to\ell_1$ which has the same distortion to $f$.
\end{prop}
By applying proposition \ref{prop:shrink}, we only have to consider when all the paths of $G$ between two endpoints of $G$ are equal.

Let $k,\ell\in\N$. We define the graph $K_{2,2k}^{\ell}$ as follows. Start with the complete bipartite graph $K_{2,2k}$ and fix its bipartition $A \cup B$ with $|A|=2$, say $A=\{0,1\}$ and $B=\{b_1,\dots,b_{2k}\}$. For each $b_i \in B$, subdivide both edges $(0,b_i)$ and $(1,b_i)$ exactly $\ell-1$ times. Equivalently, each edge $(0,b_i)$ and $(1,b_i)$ is replaced by a path of length $\ell$. In this way, every $b_i$ lies on a unique path of length $2\ell$ between $0$ and $1$, and altogether there are $2k$ internally vertex-disjoint paths of length $2\ell$ connecting $0$ and $1$. We denote the resulting graph by $K_{2,2k}^{\ell}$.

We number the vertices as follows. First, we set $A=\{0,1\}$. For each $0 \le i < 2k$, consider the corresponding path of length $2\ell$ from $0$ to $1$.   The internal vertices on this path are labeled sequentially from the side of $0$ as
\[(2\ell-1)i + j + 2, \qquad 0 \le j < 2\ell-1.\]
Thus, each path consists of the vertices
\[0 - \big((2\ell-1)i+2\big) - \big((2\ell-1)i+3\big) - \cdots - \big((2\ell-1)(i+1)+1\big) - 1,\]
and these $2k$ paths are pairwise internally disjoint.

Here is the example when $k=2,\ell=3$.
\begin{figure}[h!]
    \centering
    \begin{tikzpicture}[scale=1, every node/.style={circle, draw, inner sep=2pt}]
      \node (0) at (0,0) {0};
      \node (1) at (12,0) {1};
    
      \node (2) at (2,3) {2};
      \node (3) at (4,3) {3};
      \node (4) at (6,3) {4};
      \node (5) at (8,3) {5};
      \node (6) at (10,3) {6};
      \draw (0)--(2)--(3)--(4)--(5)--(6)--(1);
    
      \node (7) at (2,1.5) {7};
      \node (8) at (4,1.5) {8};
      \node (9) at (6,1.5) {9};
      \node (10) at (8,1.5) {10};
      \node (11) at (10,1.5) {11};
      \draw (0)--(7)--(8)--(9)--(10)--(11)--(1);
    
      \node (12) at (2,-1.5) {12};
      \node (13) at (4,-1.5) {13};
      \node (14) at (6,-1.5) {14};
      \node (15) at (8,-1.5) {15};
      \node (16) at (10,-1.5) {16};
      \draw (0)--(12)--(13)--(14)--(15)--(16)--(1);
    
      \node (17) at (2,-3) {17};
      \node (18) at (4,-3) {18};
      \node (19) at (6,-3) {19};
      \node (20) at (8,-3) {20};
      \node (21) at (10,-3) {21};
      \draw (0)--(17)--(18)--(19)--(20)--(21)--(1);
    
    \end{tikzpicture}
    \caption{Graph $K_{2,2k}^{\ell}$ where $k=2, \ell=3$}
\end{figure}
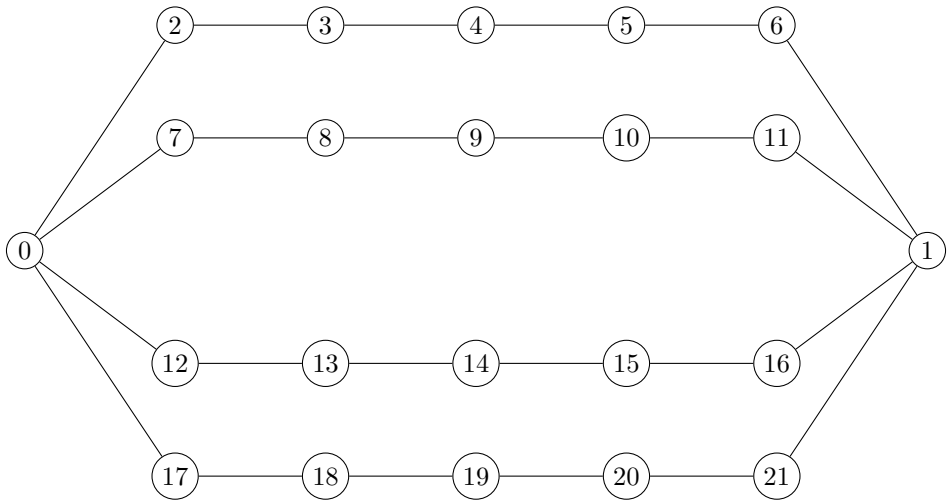

All we have to do is constructing embedding of $K_{2,2k}^\ell$ into $\ell_1$ with distortion $\frac{3k-2}{2k-1}$. To achieve this, we combine two types of cuts. Define two algorithms which generates random cuts as follows:

\begin{algorithm}[H]
  \caption{Generate Random Cut I}
  \label{alg:cut1}
  \begin{algorithmic}[1]
    \State Pick $i\in[2k]$ and $j\in[\ell]$ uniformly and independently.
    \State \textbf{return} $\{(2\ell-1)(i-1)+j+1,\ldots,(2\ell-1)(i-1)+(2\ell-j+1)\}$
  \end{algorithmic}
\end{algorithm}
\begin{algorithm}[H]
  \caption{Generate Random Cut II}
  \label{alg:cut2}
  \begin{algorithmic}[1]
    \State Pick $I\in\binom{[2k]}{k}$ and $b\in[\ell]$ uniformly and independently.
    \State Initialize $S:=\{0\}$
    \For{$i\gets 1$ to $2k$}
        \If{$i\in I$}
            \State $S:=S\cup\{(2\ell-1)(i-1)+2,\ldots,(2\ell-1)(i-1)+2\ell-b+1\}$
        \Else
            \State $S:=S\cup\{(2\ell-1)(i-1)+2,\ldots,(2\ell-1)(i-1)+b\}$
        \EndIf
    \EndFor
    \State \textbf{return} $S$
  \end{algorithmic}
\end{algorithm}

where for any integer $a,b\in\mathbb{Z}$, we write
\[\{a,\ldots,b\}=\begin{cases}
    \emptyset &(b<a)\\
    \{q\in\mathbb{Z}:a\leq q\leq b\}&(a\leq b).
\end{cases}\]

Let $S_1,S_2$ be the random set generating by our algorithm I, II respectively. For $S\subset V(K_{2,2k}^\ell)$, we write $\rho_S(x,y)=|1_S(x)-1_S(y)|$ as the cut metric on $S$. We claim the embedding into $\ell_1$ giving by in terms of cut, say
\[d_1=\frac{2\ell k(k-1)}{2k-1}\mathbb{E}[\rho_{S_1}]+2\ell\mathbb{E}[\rho_{S_2}]\]
has $\frac{3k-2}{2k-1}$ distortion.

Denote $d$ as the metric of $K_{2,2k}^\ell$. Now, for distinct $x,y\in V(K_{2,2k}^\ell)$, we evaluate $d_1(x,y)/d(x,y)$ one by one. We calculate 7 cases and the rest cases can reduce to one of the below cases by the symmetricity of our algorithms.

\newcommand{\E}{\mathbb{E}}
\noindent\textbf{Case I} $x=0,y=1$.

It is easy to verify $\E[\rho_{S_1}(x,y)]=0$ and $\E[\rho_{S_2}(x,y)]=1$, so we have $d_1(x,y)/d(x,y)=1$.

\noindent\textbf{Case II} $x=0,y\in\{2,\ldots,\ell+1\}$.

Note that $d(x,y)=y-1$. We have
\begin{align*}
    \E[\rho_{S_1}(x,y)]&=\Pr[y\in S_1]\\
    &=\frac{y-1}{2k\ell}
\end{align*}
and
\begin{align*}
    \E[\rho_{S_2}(x,y)]&=\Pr[y\notin S_2]\\
    &=\Pr[1\in I]\Pr[y\notin S_2|1\in I]+\Pr[1\notin I]\Pr[y\notin S_2|1\notin I]\\
    &=\frac{1}{2}\cdot 0+\frac{1}{2}\cdot\frac{y-1}{\ell}=\frac{y-1}{2\ell}.
\end{align*}
Thus, we have $d_1(x,y)/d(x,y)=\frac{3k-2}{2k-1}$.

\noindent\textbf{Case III} $x=0, y\in \{\ell+2,\ldots,2\ell\}$.

Clearly, $d(x,y)=y-1$. On the otherhand, we can compute
\begin{align*}
    \E[\rho_{S_1}(x,y)]&=\Pr[y\in S_1]\\
    &=\frac{2\ell-y+1}{2k\ell}
\end{align*}
and
\begin{align*}
    \E[\rho_{S_2}(x,y)]&=\Pr[y\notin S_2]\\
    &=\Pr[1\in I]\Pr[y\notin S_2|1\in I]+\Pr[1\notin I]\Pr[y\notin S_2|1\notin I]\\
    &=\frac{1}{2}\cdot 1+\frac{1}{2}\cdot\frac{y-1-\ell}{\ell}=\frac{y-1}{2\ell}.
\end{align*}
So,
\[d(x,y)=\frac{2\ell k(k-1)}{2k-1}\cdot \frac{2\ell-y+1}{2k\ell}+2\ell\cdot\frac{y-1}{2\ell}=y-1+\frac{k-1}{2k-1}(2\ell-y+1).\]
From the assumption, $0\leq 2\ell-y+1\leq y-1$ holds. So we have $1\leq d_1(x,y)/d(x,y)\leq\frac{3k-2}{2k-1}$.

\noindent\textbf{Case IV} $x,y\in \{2,\ldots,\ell+1\}$.

Without loss of generality, assume $x<y$. Then, $d(x,y)=y-x$. On the otherhand,
\begin{align*}
    \E[\rho_{S_1}(x,y)]&=\Pr[x\notin S_1,y\in S_1]\\
    &=\frac{y-x}{2k\ell}
\end{align*}
and since $y\in S_2$ always means $x\in S_2$,
\begin{align*}
    \E[\rho_{S_2}(x,y)]&=\Pr[x\in S_2,y\notin S_2]\\
    &=\Pr[1\in I]\Pr[x\in S_2,y\notin S_2|1\in I]+\Pr[1\notin I]\Pr[x\in S_2,y\notin S_2|1\notin I]\\
    &=\frac{1}{2}\cdot 0+\frac{1}{2}\cdot\frac{y-x}{\ell}=\frac{y-x}{2\ell}.
\end{align*}
we have $d_1(x,y)=\frac{3k-2}{2k-1}(y-x)$. Thus, $d_1(x,y)/d(x,y)=\frac{3k-2}{2k-1}$.

\noindent\textbf{Case V} $x\in\{2,\ldots,\ell+1\},y\in\{\ell+2,\ldots,2\ell\}$.

Note $d(x,y)=y-x$. We can calculate
\begin{align*}
    \E[\rho_{S_1}(x,y)]&=\Pr[x\notin S_1,y\in S_1]+\Pr[x\in S_1,y\notin S_1]\\
    &=\frac{|2\ell+1-x-y|}{2k\ell}
\end{align*}
and since $y\in S_2$ always means $x\in S_2$,
\begin{align*}
    \E[\rho_{S_2}(x,y)]&=\Pr[x\in S_2,y\notin S_2]\\
    &=\Pr[1\in I]\Pr[x\in S_2,y\notin S_2|1\in I]+\Pr[1\notin I]\Pr[x\in S_2,y\notin S_2|1\notin I]\\
    &=\frac{1}{2}\cdot\frac{y-\ell-1}{\ell}+\frac{1}{2}\cdot\frac{\ell-x+1}{\ell}=\frac{y-x}{2\ell}.
\end{align*}
However, $|2\ell+1-x-y|\leq y-x$ always hold, so $1\leq d_1(x,y)/d(x,y)\leq\frac{3k-2}{2k-1}$.

\noindent\textbf{Case VI} $x\in\{2,\ldots,\ell+1\},y\in \{2\ell+1,\ldots,3\ell\}$.

Clearly, $d(x,y)=(x-1)+(y-2\ell)=x+y-(2\ell+1)$. Since $x\in S_1$ always means $y\notin S_1$ and vice versa, we have
\begin{align*}
    \E[\rho_{S_1}(x,y)]&=\Pr[x\notin S_1,y\in S_1]+\Pr[x\in S_1,y\notin S_1]\\
    &=\frac{y-2\ell}{2k\ell}+\frac{x-1}{2k\ell}=\frac{(x+y)-(2\ell+1)}{2k\ell}
\end{align*}
and
\begin{align*}
    \E[\rho_{S_2}(x,y)]&=\Pr[x\notin S_2,y\in S_2]+\Pr[x\in S_2,y\notin S_2].
\end{align*}
We first compute $\Pr[x\notin S_2,y\in S_2]$.

\begin{align*}
    \Pr[x\notin S_2,y\in S_2]&=\Pr[1,2\in I]\Pr[x\notin S_2,y\in S_2|1,2\in I]\\
    &+\Pr[1\in I,2\notin I]\Pr[x\notin S_2,y\in S_2|1\in I,2\notin I]\\
    &+\Pr[1\notin I,2\in I]\Pr[x\notin S_2,y\in S_2|1\notin I,2\in I]\\
    &+\Pr[1\notin I,2\notin I]\Pr[x\notin S_2,y\in S_2|1\notin I,2\notin I]\\
    &=\frac{k}{2(2k-1)}\Pr[x\notin S_2,y\in S_2|1\notin I,2\in I]\\
    &+\frac{k-1}{2(2k-1)}\Pr[x\notin S_2,y\in S_2|1\notin I,2\notin I]\\
    &=\frac{k}{2(2k-1)}\frac{x-1}{\ell}+\frac{k-1}{2(2k-1)}\max\left\{\frac{x-y+2\ell-1}{\ell},0\right\}
\end{align*}
Similary, we have
\begin{align*}
    \Pr[x\in S_2,y\notin S_2]&=\frac{k}{2(2k-1)}\frac{y-2\ell}{\ell}+\frac{k-1}{2(2k-1)}\max\left\{\frac{y-x-2\ell+1}{\ell},0\right\}.
\end{align*}
So,
\begin{align*}
    d_1(x,y)&=\frac{2\ell k(k-1)}{2k-1}\cdot\frac{(x+y)-(2\ell+1)}{2k\ell}+2\ell\cdot\frac{k}{2(2k-1)}\cdot\frac{(x+y)-(2\ell+1)}{\ell}\\
    &+2\ell\cdot\frac{k-1}{2(2k-1)}\cdot \left|\frac{y-x-2\ell+1}{\ell}\right|\\
    &=d(x,y)+\frac{k-1}{2k-1}|y-x-2\ell+1|
\end{align*}
Since $|y-x-2\ell+1|\leq d(x,y)$, we have $1\leq d_1(x,y)/d(x,y)\leq\frac{3k-2}{2k-1}$.

\noindent\textbf{Case VII} $x\in\{2,\ldots,\ell+1\},y\in \{3\ell+1,\ldots,4\ell-1\}$.

Firstly, we have
\begin{align*}
    d(x,y)&=\min\{(x+y)-(2\ell+1),(6\ell+1)-(x+y)\}    \\
    &=2\ell-|(4\ell-y)-(x-1)|.
\end{align*}
On the other hand,
\begin{align*}
    \E[\rho_{S_1}(x,y)]&=\Pr[x\notin S_1,y\in S_1]+\Pr[x\in S_1,y\notin S_1]\\
    &=\frac{4\ell-y}{2k\ell}+\frac{x-1}{2k\ell}=\frac{(4\ell-y)+(x-1)}{2k\ell}
\end{align*}
and
\begin{align*}
    \E[\rho_{S_2}(x,y)]&=\Pr[x\notin S_2,y\in S_2]+\Pr[x\in S_2,y\notin S_2].
\end{align*}
Let us evaluate each term one by one. First,
\begin{align*}
    \Pr[x\notin S_2,y\in S_2]&=\Pr[1,2\in I]\Pr[x\notin S_2,y\in S_2|1,2\in I]\\
    &+\Pr[1\in I,2\notin I]\Pr[x\notin S_2,y\in S_2|1\in I,2\notin I]\\
    &+\Pr[1\notin I,2\in I]\Pr[x\notin S_2,y\in S_2|1\notin I,2\in I]\\
    &+\Pr[1\notin I,2\notin I]\Pr[x\notin S_2,y\in S_2|1\notin I,2\notin I]\\
    &=\frac{k}{2(2k-1)}\cdot\Pr[x\notin S_2,y\in S_2|1\notin I,2\in I]\\
    &=\frac{k}{2(2k-1)}\Pr[b\in [1, x-1]\cap [1,4\ell-y]|1\notin I,2\in I]\\
    &=\frac{k}{2(2k-1)}\cdot\frac{\min\{x-1, 4\ell-y\}}{\ell}
\end{align*}
and second,
\begin{align*}
    \Pr[x\in S_2,y\notin S_2]&=\Pr[1,2\in I]\Pr[x\in S_2,y\notin S_2|1,2\in I]\\
    &+\Pr[1\in I,2\notin I]\Pr[x\in S_2,y\notin S_2|1\in I,2\notin I]\\
    &+\Pr[1\notin I,2\in I]\Pr[x\in S_2,y\notin S_2|1\notin I,2\in I]\\
    &+\Pr[1\notin I,2\notin I]\Pr[x\in S_2,y\notin S_2|1\notin I,2\notin I]\\
    &=\frac{k-1}{2(2k-1)}\Pr[b\in[4\ell-y+1,\ell]]+\frac{k}{2(2k-1)}\cdot 1\\
    &+\frac{k}{2(2k-1)}\Pr[b\in [x,\ell]\cap [4\ell-y+1,\ell]]+\frac{k-1}{2(2k-1)}\Pr[b\in [x,\ell]]\\
    &=\frac{k-1}{2(2k-1)}\frac{y-3\ell}{\ell}+\frac{k}{2(2k-1)}+\frac{k-1}{2(2k-1)}\frac{\ell-x+1}{\ell}\\
    &+\frac{k}{2(2k-1)}\frac{\ell-\max\{4\ell-y,x-1\}}{\ell}\\
    &=\frac{k-1}{2(2k-1)}\frac{(y-2\ell)-(x-1)}{\ell}+\frac{k}{2k-1}\\
    &-\frac{k}{2(2k-1)}\frac{\max\{4\ell-y,x-1\}}{\ell}
\end{align*}
Thus, we have
\begin{align*}
    \E[\rho_{S_2}(x,y)]=\frac{k-1}{2(2k-1)}\frac{(y-2\ell)-(x-1)}{\ell}+\frac{k}{2k-1}-\frac{k}{2(2k-1)}\frac{|(4\ell-y)-(x-1)|}{\ell}
\end{align*}
Therefore,
\begin{align*}
    d(x,y)&=\frac{k-1}{2k-1}((4\ell-y)+(x-1)|-\frac{k}{2k-1}|(4\ell-y)-(x-1)\\
    &+\frac{2\ell k}{2k-1}+\frac{k-1}{2k-1}((y-2\ell)-(x-1))\\
    &=2\ell-\frac{k}{2k-1}|(4\ell-y)-(x-1)|
\end{align*}
Since $|(4\ell-y)-(x-1)|\leq \ell$, one can check $1\leq d(x,y)/d_1(x,y)\leq \frac{3k-2}{2k-1}$.
\bibliographystyle{plain}
\bibliography{ref}

@inproceedings{Chakrabarti_et_al_2008,
author = {Chakrabarti, Amit and Jaffe, Alexander and Lee, James R. and Vincent, Justin},
title = {Embeddings of Topological Graphs: Lossy Invariants, Linearization, and 2-Sums},
year = {2008},
isbn = {9780769534367},
publisher = {IEEE Computer Society},
address = {USA},
url = {https://doi.org/10.1109/FOCS.2008.79},
doi = {10.1109/FOCS.2008.79},
booktitle = {Proceedings of the 2008 49th Annual IEEE Symposium on Foundations of Computer Science},
pages = {761–770},
numpages = {10},
series = {FOCS '08}
}

@Article{Gupta2004,
author={Gupta, Anupam
and Newman, Ilan
and Rabinovich, Yuri
and Sinclair, Alistair},
title={Cuts, Trees and  $\ell_1$-Embeddings of Graphs},
journal={Combinatorica},
year={2004},
month={Apr},
day={01},
volume={24},
number={2},
pages={233-269},
issn={1439-6912},
doi={10.1007/s00493-004-0015-x},
url={https://doi.org/10.1007/s00493-004-0015-x}
}

@article{ostrovsky2007low,
  title={Low distortion embeddings for edit distance},
  author={Ostrovsky, Rafail and Rabani, Yuval},
  journal={Journal of the ACM (JACM)},
  volume={54},
  number={5},
  pages={23--es},
  year={2007},
  publisher={ACM New York, NY, USA}
}

@book{Deza,
  title={Geometry of cuts and metrics},
  author={Michel Deza and Monique Laurent},
  booktitle={Algorithms and Combinatorics},
  year={1997},
  publisher={springer}
}
\end{document}